\documentclass[12pt]{amsart}
\usepackage{amscd}      
\usepackage{amssymb}
\usepackage{amsmath, amsthm, graphics}
\usepackage{xypic}      
\LaTeXdiagrams          
\usepackage[all]{xy}
\xyoption{2cell} \UseAllTwocells \xyoption{frame} \CompileMatrices
\allowdisplaybreaks[3]

\usepackage{latexsym}
\usepackage{epsfig}
\usepackage{amsfonts}
\usepackage{enumerate}
\usepackage{times}

\newtheorem{prop}{Proposition}[section]

\newtheorem{question}[prop]{Question}
\newtheorem{conjecture}[prop]{Conjecture}

\newtheorem{theorem}{Theorem}[section]

\theoremstyle{remark}

\theoremstyle{remark}

\numberwithin{equation}{section}

\newcommand{\Mbar}{\overline{\M}}

\newcommand{\X}{\mathcal{X}}
\newcommand{\Y}{\mathcal{Y}}

\newcommand{\M}{\mathcal{M}}
\newcommand{\C}{\mathcal{C}}

\newcommand{\B}{\mathcal{B}}

\newcommand{\D}{\mathcal{D}}

\newcommand{\sL}{\mathcal{L}}

\def\<{\left\langle}
\def\>{\right\rangle}

\begin{document}

\title{On the geometry of orbifold Gromov-Witten invariants}

\author[Tseng]{Hsian-Hua Tseng}
\address{Department of Mathematics\\ Ohio State University\\ 100 Math Tower, 231 West 18th Ave. \\ Columbus \\ OH 43210\\ USA}
\email{hhtseng@math.ohio-state.edu}

\date{\today}

\begin{abstract}
We consider the question of how geometric structures of a Deligne-Mumford stack affect its Gromov-Witten invariants. The two geometric structures studied here are {\em gerbes} and {\em root constructions}. In both cases, we explain conjectures on Gromov-Witten theory for these stacks and survey some recent progress on these conjectures.
\end{abstract}

\maketitle

\section{Introduction}
\subsection{Overview}
Gromov-Witten theory of Deligne-Mumford stacks was introduced in \cite{agv1, agv2} after the work \cite{cr} in symplectic setting. It has been an active area of research in the past decade, in particular in connection with the {\em crepant resolution conjecture} \cite{cit}, \cite{c_r}, \cite{b_g} and with mirror symmetry (see e.g. \cite{ccit}). 

In \cite{tseng_iccm10}, the author gave a survey on some aspects of orbifold Gromov-Witten theory. This paper can be considered as a sequel to \cite{tseng_iccm10}: here we further elaborate on some questions discussed in \cite{tseng_iccm10} and give an update on recent progress. In Section \ref{subsec:geom} we explain some natural geometric structures of Deligne-Mumford stacks. In Section \ref{subsec:str}, we indicate that these structures lead to the following two questions: 
\begin{enumerate}
\item
 the structure of Gromov-Witten theory of gerbes, which we discuss in Section \ref{sec:gerbe}. 
\item
the determination of the Gromov-Witten theory of a root stack, which we discuss in Section \ref{sec:roots}.
\end{enumerate}
Basic definitions and constructions of orbifold Gromov-Witten theory are given in many places, for example \cite{agv1, agv2}, \cite{cr}, and \cite{tseng_iccm10}, and are reviewed in Section \ref{sec:GW} below. Throughout this paper, we work over $\mathbb{C}$. 

\subsection{Geometry of Deligne-Mumford stacks}\label{subsec:geom}
We discuss the geometric structure of a smooth Deligne-Mumford stack in relation to its coarse moduli space.

Let $\X$ be a smooth Deligne-Mumford stack. By the main result of \cite{k_m}, $\X$ admits a coarse moduli space $X$ which we assume to be a scheme. The canonical morphism $$\X\to X$$
is universal in the sense that any morphism $\X\to Z$ with $Z$ a scheme must factor through $X$, i.e. there exists a morphism $X\to Z$ such that the following diagram is commutative:
\begin{displaymath}
    \xymatrix{ 
    \X \ar[dr]\ar[d] & \\
   X\ar[r]& Z.}
\end{displaymath}
Let $G$ be the generic stabilizer group of $\X$. The construction of {\em rigidification} in \cite{acv} allows one to ``remove'' the group $G$ from stabilizers. More precisely, this yields a stack $\X_1$ together with a morphism\footnote{This morphism is also constructed in \cite{bn}.} $$\X\to \X_1$$ such that the generic stabilizer group of $\X_1$ is $G/G$, i.e. the trivial group, and $\X_1$ also has $X$ as its coarse moduli space. We thus have a factorization
\begin{displaymath}
    \xymatrix{ 
    \X \ar[dr]\ar[d] & \\
  \X_1\ar[r]& X.}
\end{displaymath}
By construction, the morphism $\X\to \X_1$ has the structure of a {\em gerbe}. The stacky locus (i.e. the locus where the stabilizer groups are not trivial) of $\X_1$ is of codimension at least $1$. 

By assumption, the scheme $X$ has only quotient singularities. \'Etale locally $X$ is a quotient of a smooth scheme by a finite group. These local presentations of $X$ assemble to an \'etale groupoid. Let $$\X_{can}$$ denote the stack quotient of this groupoid. The coarse moduli space of $\X_{can}$ is also $X$, and we have a canonical morphism $$\X_{can}\to X.$$ $\X_{can}$ is the {\em canonical stack} associated to $X$, in the sense of \cite[Section 4.1]{fmn}. By \cite[Theorem 4.6]{fmn}, for any Deligne-Mumford stack $\Y$ with trivial generic stabilizer and a dominant codimension preserving morphism $\Y\to X$, there exists a unique factorization 
\begin{displaymath}
    \xymatrix{ 
    \Y \ar[r]\ar[d] &\X_{can}\ar[dl] \\
  X.&}
\end{displaymath}
In particular, there is a morphism $$\X_1\to \X_{can}.$$ Therefore we have a factorization of $\X\to X$:
\begin{equation}\label{eqn:factor}
\X\to \X_1\to \X_{can}\to X.
\end{equation}
By \cite{gs}, under some mild hypothesis, the map $\X_1\to \X_{can}$ is a composition of the {\em root construction} of \cite{agv2} and \cite{ca}. A review of root constructions is given in Section \ref{subsec:conj_codim1}.

\subsection{Structures of Gromov-Witten theory}\label{subsec:str}

In this paper, we take the following perspective: the factorization (\ref{eqn:factor}) leads to structural results of Gromov-Witten invariants of $\X$. 

\subsubsection{Stage $1$}
The first stage of the factorization, $\X\to \X_1$, is a gerbe. In Section \ref{sec:gerbe} below we explain the conjectural description of Gromov-Witten theory of $\X$ obtained from this gerbe structure, see Conjecture \ref{conj:gerbe}. As evidence supporting the conjecture, we also explain in some details the state-of-the-art result on the structure of Gromov-Witten invariants of a large class of gerbes whose bands are trivial.

\subsubsection{Stage $2$}
The second stage of the factorization, $\X_1\to \X_{can}$, is a composition of root constructions. It is natural\footnote{Apparently experts have believed this to be true.} to expect that the Gromov-Witten theory of a root construction can be explicitly determined by the input data, namely the base stack and the divisor along which the root is taken. However, such a result is not proven in full generality. In Section \ref{sec:roots}, we formulate this expectation more precisely as Conjecture \ref{conj:roots}. We also discuss some details of a proof of Conjecture \ref{conj:roots} in a special case\footnote{Despite the expectation, even this special case was proven very recently, in 2015.}.

\subsubsection{Stage $3$}
The third stage of this factorization, $\X_{can}\to X$, is canonical: as mentioned above, the stack $\X_{can}$ is canonically associated to the scheme $X$. By construction the stacky locus of $\X_{can}$ is of codimension at least $2$. Therefore, the morphism $\X_{can}\to X$ is birational and {\em crepant} (i.e. it identifies the canonical divisors $K_{\X_{can}}$ and $K_X$ via pull-back). As advocated by M. Satriano and others, $\X_{can}\to X$ should be considered as a {\em stacky (crepant) resolution}. The so-called {\em crepant resolution conjecture}, originated from String Theory, predicts an equivalence between the Gromov-Witten theory of $\X_{can}$ and Gromov-Witten theory of other crepant resolutions of $X$ (if existed). As noted above, crepant resolution conjecture has been the main driving force for much of the developments in orbifold Gromov-Witten theory in the past decade. In this paper we do not attempt to give a survey on the status of crepant resolution conjecture. The papers \cite{c_r} and \cite{cij} contains nice discussions about this conjecture, to which we refer the readers.

\subsection{Acknowledgment}
The author is very grateful to Xiang Tang and Fenglong You for inspiring and fruitful collaborations on results surveyed in this paper. The author is supported in part by NSF grant DMS-1506551 and a Simons Foundation Collaboration Grant.

\section{Gromov-Witten theory of Deligne-Mumford stacks}\label{sec:GW}
In this Section we give a review on the basic constructions of Gromov-Witten theory of Deligne-Mumford stacks. In addition to the original papers \cite{cr}, \cite{agv1}, and \cite{agv2} mentioned above, expository accounts of Gromov-Witten theory of Deligne-Mumford stacks can be found in \cite{ab} and \cite{tseng05}.

Let $\X$ be a smooth proper Deligne-Mumford stack over $\mathbb{C}$ with projective coarse moduli space $X$. The Gromov-Witten theory of $\X$ is based on moduli stacks of stable \footnote{These are often called {\em twisted stable maps} \cite{av} or {\em orbifold stable maps} \cite{cr}.} to $\X$. Let $g, n\in \mathbb{Z}_{\geq 0}$ and let $d\in H_2(X, \mathbb{Q})$. The moduli stack $$\Mbar_{g,n}(\X, d)$$
of $n$-pointed genus $g$ stable maps to $\X$ of degree $d$ parametrizes objects of the following kind: $$(f:\C\to \X, p_1,...,p_n)$$ where

\begin{enumerate}
\item
$\C$ is a nodal orbifold curve. This means that $\C$ is a connected nodal curve with isolated stacky points of the following two types:
\begin{enumerate}
\item
If a stacky point is not nodal, then it has an \'etale neighborhood of the form $$[\text{Spec}\mathbb{C}[x]/\mu_r]$$ for some $\mu_r$ where $\zeta\in \mu_r$ acts on the coordinate $x$ by multiplication.

\item
If a stacky point is nodal, then it has an \'etale neighborhood of the form $$[(\text{Spec}\mathbb{C}[x,y]/(xy))/\mu_r]$$ for some $\mu_r$ where $\zeta\in \mu_r$ acts on the coordinates\footnote{The action of $\mu_r$ on coordinates $x, y$ considered here, which leaves $xy$ invariant, is called {\em balanced}. Such a node is called a balanced node. The more general notion of stable maps to $\X$ introduced in \cite{av} allows nodes which are not balanced. However, for Gromov-Witten theory we require the nodes to be balanced.} $x, y$ via $$x\mapsto \zeta x, \quad y\mapsto \zeta^{-1} y.$$
\end{enumerate}

\item
$f:\C\to \X$ is a representable morphism of stacks. This means that for $x\in \C$, the group homomorphism between stabilizer groups induced by $f$, $$Aut(x)\to Aut(f(x)),$$
is injective.

\item
$p_1,..., p_n\in \C$ are distinct (possibly stacky) points which are not nodal.

\item
$f_*[\C]=d$. Here we use the fact that $H_2(\X, \mathbb{Q})\simeq H_2(X, \mathbb{Q})$.

\item
The data $(f:\C\to \X, p_1,...,p_n)$ have only finitely many automorphisms. Combinatorially this means the following:
\begin{enumerate}
\item
On an irreducible component of $\C$ of genus $0$, the number of special points (i.e. marked points and nodes) is at least $3$.
\item
On an irreducible component of $\C$ of genus $1$, the number of special points is at least $1$.

\end{enumerate}

\end{enumerate}

Furthermore, we require a technical condition: for any scheme $S$, on an $S$-family of stable maps to $\X$, 
\begin{displaymath}
    \xymatrix{ 
    \C \ar[r]\ar[d]& \X\\
  S,& }
\end{displaymath}
the marked points of $\C$ form gerbes over $S$ which we require to have sections. This condition, used in \cite{agv1} and \cite{tseng05}, is related to the definition of evaluation maps, which we discuss below.

By the main results of \cite{av}, the moduli stack $\Mbar_{g,n}(\X,d)$ is a proper Deligne-Mumford stack with projective coarse moduli space. It admits a universal family of stable maps,
\begin{displaymath}
    \xymatrix{ 
    \C \ar[r]^{f}\ar[d]_{\pi}& \X\\
  \Mbar_{g,n}(\X,d).& }
\end{displaymath}
In general, the stack $\Mbar_{g,n}(\X,d)$ can be pretty bad: it can be arbitrarily singular, with many components of various dimensions, etc. To obtain intersection numbers with reasonable properties, we need to use virtual fundamental classes. As in \cite{agv2}, the object $R^*\pi_*f^*T_\X$ in the derived category of $\Mbar_{g,n}(\X,d)$ yields a perfect obstruction theory in the sense of \cite{bf}. By the recipe of \cite{bf}, this yields a virtual fundamental class $$[\Mbar_{g,n}(\X,d)]^w\in H_*(\Mbar_{g,n}(\X,d),\mathbb{Q}).$$
We should point out that, because of the technical condition on marked points mentioned above, this virtual fundamental class needs a minor modification. See \cite{agv1} and \cite{tseng05} for details. 

We now discuss evaluation maps. Given a stable map $f:\C\to \X$, restricting to a marked point $p_i\in \C$ yields a morphism $p_i\to \X$. Since $p_i$ is possibly a stacky point $[\text{pt}/\mu_r]$, the morphism $p_i\to \X$ gives an object of the {\em inertia stack} of $\X$, $$I\X:=\X\times_{\Delta, \X\times \X, \Delta} \X,$$
where $\Delta: \X\to \X\times \X$ is the diagonal morphism. Because of the above technical condition\footnote{Without this technical condition, the evaluation map takes values in the {\em rigidified} inertia stack, see \cite{agv1}, \cite{agv2} for details. For the purpose of this paper it is more convenient to work with inertia stacks, hence we impose this technical condition.} on marked points, this restriction construction yields a morphism $$ev_i: \Mbar_{g,n}(\X, d)\to I\X,$$ called the $i$-th evaluation map. 

The $i$-th tautological line bundle on $\Mbar_{g,n}(\X,d)$, 
$$L_i\to \Mbar_{g,n}(\X,d),$$
has fibers over a stable map $(f:\C\to \X, p_1,...,p_n)$ the cotangent line $T^*_{\bar{p}_i}C$ of the {\em coarse curve} $C$ at the (image of) $i$-th marked point. The descendant classes on $\Mbar_{g,n}(\X,d)$ are first Chern classes of $L_i$, $$\psi_i:=c_1(L_i), \quad 1\leq i\leq n.$$

We are now ready to define Gromov-Witten invariants of $\X$. Given integers $k_1, ..., k_n\geq 0$ and classes in the {\em Chen-Ruan orbifold cohomology of $\X$}, $a_1,..., a_n\in H^*(I\X)$, we define the following descendant Gromov-Witten invariant of $\X$, 
\begin{equation}
\<\prod_{i=1}^n a_i\psi_i^{k_i}\>_{g,n,\beta}^\X:=\int_{[\Mbar_{g,n}(\X,d)]^w} \prod_{i=1}^n ev_i^*(a_i)\psi_i^{k_i}\in \mathbb{Q}.
\end{equation}

\section{Codimension $0$: gerbes}\label{sec:gerbe}
In this Section we discuss Gromov-Witten theory of gerbes. Our study is guided by a conjecture originated in Physics \cite{hhps} and further developed in \cite{tt10}. This is Conjecture \ref{conj:gerbe} below.

\subsection{The conjecture}
Essentially a gerbe is a stack with non-trivial generic stabilizers. Let $\Y$ be a smooth Deligne-Mumford stack whose generic stabilizer is a finite group $G$. Then there is a map $$\pi:\Y\to \B$$ where $\B$ is a smooth Deligne-Mumford stack with trivial generic stabilizer. The map $\pi$ is known as a {\em $G$-gerbe} over $\B$. The only $G$-gerbe over $\B=\text{pt}$ is $$BG=[\text{pt}/G].$$ Given an open cover $\{U_i\}$ of $\B$, one can construct $\Y$, the total space of a $G$-gerbe over $\B$, by gluing $U_i\times BG$ using the following data:
\begin{equation}
\begin{split}
&\phi_{ij}\in Aut(G) \quad \text{for each double overlap } U_{ij}:=U_i\cap U_j,\\
&g_{ijk}\in G \quad \text{for each triple overlap } U_{ijk}:=U_i\cap U_j\cap U_k.
\end{split}
\end{equation}
The following requirements are imposed:
\begin{equation}\label{eqn:gerbe_camp}
\begin{split}
&\phi_{jk}\circ \phi_{ij}=\text{Ad}_{g_{ijk}}\circ \phi_{ik} \quad \text{on } U_{ijk}\\
&g_{jkl}g_{ijl}=\phi_{kl}(g_{ijk})g_{ikl} \quad \text{on } u_i\cap U_j\cap U_k\cap U_l.
\end{split}
\end{equation}
Let $Out(G):=Aut(G)/Inn(G)$ be the group of outer automorphisms of $G$. By (\ref{eqn:gerbe_camp}), the images of $\phi_{ij}$ under the quotient map $Aut(G)\to Out(G)$ satisfy the usual compatibility requirement for principal bundles and can be used to construct a principal $Out(G)$-bundle over $\B$, which we denote by $$\overline{\Y}\to \B.$$ This principal bundle $\overline{\Y}\to \B$ is called the {\em band} of the $G$-gerbe $\pi$. As in \cite{hhps}, for a $G$-gerbe $\pi: \Y\to \B$, one can construct a pair $$(\widehat{\Y}, c),$$ where $\widehat{\Y}$ is a disconnected space with an \'etale map $\widehat{\pi}:\widehat{\Y}\to \B$ and $c$ is a $U(1)$-valued $2$-cocycle on $\widehat{\Y}$, see also \cite{tt10}. More precisely, $$\widehat{\Y}:=[\overline{\Y}\times \widehat{G}/Out(G)].$$ Here $\widehat{G}$ is the set of isomorphism classes of irreducible complex representations of $G$, treated as a space of isolated points. The group $Out(G)$ acts on $\widehat{G}$ by precompositions and acts on $\overline{\Y}\times \widehat{G}$ diagonally. The $U(1)$-valued $2$-cocycle $c$ measures the failure of certain family of vector spaces to be a locally free sheaf on $\widehat{\Y}$, see \cite{tt10} for more details. The Gromov-Witten theory of $\widehat{\Y}$ incorporating a twist by the cocycle $c$ can be obtained from the general construction of \cite{pry}. 

\begin{conjecture}[see \cite{hhps} and \cite{tt10}]\label{conj:gerbe}
For a $G$-gerbe $\Y\to \B$, the Gromov-Witten theory of $\Y$ is equivalent to the Gromov-Witten theory of $\widehat{\Y}$ twisted by $c$.
\end{conjecture}
Upon first look, Conjecture \ref{conj:gerbe} seems strange and unmotivated. It is thus surprising that Conjecture \ref{conj:gerbe} turns out to be correct in all examples studied so far: it is proven for specific classes of examples in \cite{j}, \cite{ajt1, ajt2, ajt3}, and \cite{tt10, tt16}. In \cite{tt16} Conjecture \ref{conj:gerbe} is proven for $G$-gerbes $\Y\to \B$ whose bands $\overline{\Y}\to \B$ are trivializable principal bundles. This is so far the largest class of examples. In the rest of this Section we give a survey on the results of \cite{tt16}. 

\subsection{The trivial band case}
When the band $\overline{\Y}\to \B$ is trivializable, i.e. $\overline{\Y}\simeq \B\times Out(G)$, the construction implies that $$\widehat{\Y}\simeq \B\times \widehat{G}.$$ In other words, $\widehat{\Y}$ is a disjoint union of $|\widehat{G}|$ copies of $\B$. Furthermore, on each copy of $\B$ in $\widehat{\Y}$, the map $\widehat{\pi}$ restricts to the identity map on $\B$. For $[\rho]\in \widehat{G}$, we write $$\B_\rho\subset \widehat{\Y}$$ for the copy of $\B$ indexed by $[\rho]$. The restriction of the cocycle $c$ onto $\B_\rho$, which we denote by $c_\rho$, defines a cohomology class $$[c_\rho]\in H^2(\B, U(1)).$$ This class has a very simple description as follows. By the classification of gerbes with trivial bands (see e.g. \cite{g}), the $G$-gerbe $\pi: \Y\to \B$ determines a class $$[\pi]\in H^2(\B, Z(G)),$$ where $Z(G)\subset G$ is the center of $G$. For an irreducible representation $\rho: G\to GL(V_\rho)$, by Schur's lemma, the restriction $\rho|_{Z(G)}: Z(G)\to GL(V_\rho)$ determines a homomorphism $\rho|_{Z(G)}: Z(G)\to U(1)$. Using this to change coefficients, we obtain a map $$H^2(\B, Z(G))\to H^2(\B, U(1)).$$ The class $[c_\rho]$ is the image of $[\pi]$ under this map. 

The Gromov-Witten theory of $\widehat{\Y}$ twisted by $c$ is understood to be the Gromov-Witten theory of $\B_\rho$ twisted by $c_\rho$ for all $[\rho]\in \widehat{G}$ put together. For the precise meaning of this, see Theorem \ref{thm:gerbe} below. By the construction of \cite{pry}, the Gromov-Witten theory of $\B_\rho$ twisted by $c_\rho$ depends only on the class $[c_\rho]$, not the cocycle $c_\rho$. 

The ``state space'' of the Gromov-Witten theory of $\Y$ is the Chen-Ruan orbifold cohomology $H^*_{CR}(\Y)$. Additively, $$H^*_{CR}(\Y):=H^*(I\Y, \mathbb{C})$$ is the cohomology of $I\Y:=\Y\times_{\Y\times \Y} \Y$, the inertial stack of $\Y$. The state space of the Gromov-Witten theory of $\B$ twisted by $c_\rho$ is the twisted orbifold cohomology $H^*_{orb}(\B, c_\rho)$ of $(\B, c_\rho)$ introduced in \cite{r}. Additively, $$H^*_{orb}(\B,c_\rho):=H^*(I\B, \sL_{c_\rho})$$ is the cohomology of the inertial stack $I\B$ with coefficients in a line bundle $\sL_{c_\rho}\to I\B$ known as inner local system. As in \cite{r}, the line bundle $\sL_{c_\rho}$ can be explicitly constructed\footnote{The isomorphism class of $\sL_{c_\rho}$ admits a formal description. There is a natural map $H^2(\B, U(1))\to H^1(L\B, U(1))$, where $L\B$ is the loop space of $\B$. The inertia stack $I\B$ is contained in $L\B$ as the locus of constant loops. The image of the class of $c_\rho$ restricts to the class of $\sL_{c_\rho}$ via the inclusion $I\B\subset L\B$.} using $c_\rho$ as a cocycle. 

By a result\footnote{This isomorphism $I$ is also constructed in \cite{tt10} for $G$-gerbes whose bands are not trivial. The construction of $I$ is rather complicated, as it involves Morita equivlanece between explicit groupoid algebras. However, an explicit formula of $I$ was found in \cite{tt10}.} of \cite{tt10}, there is an isomorphism of graded vector spaces,
$$I: H^*_{CR}(\Y)\to \oplus_{[\rho]\in \widehat{G}}H^*_{orb}(\B, c_\rho).$$
For $\delta\in H^*_{CR}(\Y)$, we write $I(\delta)=\sum_{[\rho]\in \widehat{G}}I(\delta)_\rho$ where $I(\delta)_\rho\in H^*_{orb}(\B, c_{\rho})$. 

The main result of \cite{tt16} is the following, which confirms Conjecture \ref{conj:gerbe} for $G$-gerbes with trivial bands in a precise form. 

\begin{theorem}[\cite{tt16}, Theorem 1.1]\label{thm:gerbe}
Let $\Y\to \B$ be a $G$-gerbe over $\B$ whose band is trivializable. Let $g\geq 0$, $n>0$, and $a_1,...,a_n\geq 0$ be integers. Let $\beta\in H_2(\B, \mathbb{Q})$ be a curve class. Let $\delta_1,...,\delta_n\in H^*_{CR}(\Y)$. Then the following equality of descendant Gromov-Witten invariants holds:
\begin{equation}
\<\prod_{j=1}^n\delta_j \psi_j^{a_j} \>_{g,n,\beta}^\Y=\sum_{[\rho]\in \widehat{G}}\left(\frac{\text{dim} V_\rho}{|G|} \right)^{2-2g}\<\prod_{j=1}^nI(\delta_j)_\rho \psi_j^{a_j}\>_{g,n,\beta}^{\B, c_\rho}.
\end{equation}
\end{theorem}
In what follows we explain some ideas involved in the proof of Theorem \ref{thm:gerbe}. Gromov-Witten invariants of $\Y$ are intersection numbers on moduli spaces $$\Mbar_{g,n}(\Y, \beta)$$ of stable maps to $\Y$. Gromov-Witten invariants of $\B$ twisted by $c_\rho$ are intersection numbers on moduli spaces $$\Mbar_{g,n}(\B, \beta)$$ of stable maps to $\B$. Given a stable map $\C\to \Y$, one can obtain a stable map to $\B$ by stabilizing the composition $\C\to \Y\overset{\pi}{\longrightarrow} \B$. This construction defines a map $$\pi_{g,n,\beta}:\Mbar_{g,n}(\Y, \beta)\to \Mbar_{g,n}(\B, \beta).$$
The first ingredient in the proof of Theorem \ref{thm:gerbe} is the {\em virtual pushforward property} for $\pi_{g,n,\beta}$. More precisely, we prove that the pushforward $$\pi_{g,n,\beta *}[\Mbar_{g,n}(\Y, \beta)]^{w}\in H_*(\Mbar_{g,n}(\B, \beta), \mathbb{Q})$$ can be written as a linear combination of components of the virtual fundamental class $[\Mbar_{g,n}(\B, \beta)]^{w}$. To compute coefficients of this linear combination, we need to study the degrees of $\pi_{g,n,\beta}$ over components of $\Mbar_{g,n}(\B, \beta)$. More precisely, given a stable map $f:\C\to \B$ with $\C$ nonsingular, we need to classify all possible commutative diagrams
\begin{displaymath}
    \xymatrix{ 
    \C' \ar[r]^{f'}\ar[d] & \Y\ar[d]^{\pi}\\
   \C\ar[r]^{f}& \B}
\end{displaymath}
where $f': \C'\to \Y$ is a stable map and $\C'\to\C$ is an isomorphism over the non-stacky locus $\C^{ns}$.  The class of the pull-back gerbe $$f^*\Y\to \C$$ lies in $H^2(\C, Z(G))$. Let $$\C^0\subset \C$$ be obtained from $\C$ by removing small disks centered at marked points. Since $\C^0$ is homotopic to its $1$-skeleton, we have $$H^2(\C^0, Z(G))=0.$$ So the restriction of $f^*\Y\to \C$ to $\C^0$ is a trivializable $G$-gerbe. This simple observation leads to a characterization of all possible $f'$ in terms of conditions on orbifold structures on $\C'$. Consequently, the count of all such $f'$ can be obtained from the Frobenius-Mednyh formula \cite{me}. Solution to this counting problem can be extended to nodal $\C$ by a simple splitting argument.

Now let $\varphi_1,...,\varphi_n\in H^*_{orb}(\B, c_\rho)$. By definition we have $$\<\prod_{j=1}^n\varphi_j\psi_j^{a_j}\>_{g,n, \beta}^{\B,c_\rho}:=\int_{[\Mbar_{g,n}(\B, \beta)]^{w}}\theta_*(\prod_{j=1}^nev_j^*(\varphi_j))\cup \prod_{j=1}^n\psi_j^{a_j},$$
here $\theta: \otimes_{j=1}^n ev_j^*\sL_{c_\rho}\to \underline{\mathbb{C}}$ is the trivialization appearing in the construction of \cite{pry}. As discussed above, there is an equality\footnote{At the risk of not being modest enough, we want to point out that the determination of this degree function $d_{g,n,\beta}$ in \cite{tt16} is quite difficult. Also, the degrees appears in previous applications of virtual pushforwards are very simple. This makes our application of virtual pushforwards technically very sophisticated.} $$\pi_{g,n,\beta *}[\Mbar_{g,n}(\Y, \beta)]^{w}=d_{g,n,\beta} [\Mbar_{g,n}(\B, \beta)]^{w}$$ for a certain degree function $d_{g,n,\beta}$. Therefore we have $$\<\prod_{j=1}^n\varphi_j\psi_j^{a_j}\>_{g,n, \beta}^{\B,c_\rho}=\frac{1}{d_{g,n,\beta}}\int_{[\Mbar_{g,n}(\Y, \beta)]^{w}}\pi_{g,n,\beta}^*\theta_*(\prod_{j=1}^nev_j^*(\varphi_j))\cup \prod_{j=1}^n\psi_j^{a_j},$$
(note that the descendant classes $\psi_j$ are compatible under various pull-backs). We need to compare this with $$\<\prod_{j=1}^nI^{-1}(\varphi_j) \psi_j^{a_j} \>_{g,n,\beta}^\Y=\int_{[\Mbar_{g,n}(\Y, \beta)]^{w}}\prod_{j=1}^nev_j^*(I^{-1}(\varphi_j))\psi_j^{a_j}.$$
When $G$ is abelian, this can be done easily from an explicit formula for the inverse map $I^{-1}$. For non-abelian $G$, the short exact sequence $$1\to Z(G)\to G\to G/Z(G)=:K\to 1$$ leads to a factorization of $\pi:\Y\to \B$ into a composition $$\Y\to \Y'\to \B$$ where $\Y\to \Y'$ is a $Z(G)$-gerbe with trivial band and $\Y\to \B$ is a $K$-gerbe with trivial band. Theorem \ref{thm:gerbe} for $\Y\to \Y'$ is proven first. Since the center of $K$ is trivial, the $K$-gerbe $\Y'\to \B$ is necessarily trivial, and Theorem \ref{thm:gerbe} for $\Y'\to \B$ follows from either the results on Gromov-Witten theory of trivial gerbes in \cite{ajt1}, or by a direct calculation using the virtual pushforward property. Theorem \ref{thm:gerbe} for $\Y\to \B$ then follows by combining these two cases.

\subsection{Outlook}
Naturally, the next goal is to prove Conjecture \ref{conj:gerbe} for $G$-gerbes with non-trivial bands. It is hopeful that the approach used to prove Theorem \ref{thm:gerbe} can be applied to this general setting as well. However, there are several difficulties:
\begin{enumerate}
\item
If the band of $\Y\to \B$ is not trivializable, then components of the ``dual'' $\widehat{\Y}$ are not necessarily isomorphic to $\B$. This further complicates the study of Gromov-Witten theory of $\widehat{\Y}$, twisted by $c$ or not.

\item
In case of non-trivial bands, the class of $c$ does not seem to have a simple description. This makes the studt of Gromov-Witten theory of $\widehat{\Y}$ twisted by $c$ much harder.

\item
If the band is not trivial, then computation of the degree of the map $\pi_{g,n,\beta}:\Mbar_{g,n}(\Y, \beta)\to \Mbar_{g,n}(\B, \beta)$ must be more involved. While the general idea used in the trivial band case still seems hopeful, it is not clear how to, for example, take into account the information about the band.

\end{enumerate}

\section{Codimension $1$: root stacks}\label{sec:roots}
In this Section, we discuss Gromov-Witten theory of root constructions.

\subsection{Conjecture}\label{subsec:conj_codim1}
A way to introduce stack structures in codimension $1$ is the construction of root stacks, introduced in \cite{ca} and \cite[Appendix B]{agv2} (see also \cite{mo}). We review the construction. Let $\X$ be a Deligne-Mumford stack and $$\D\subset \X$$ a smooth irreducible divisor. Given a positive integer $r$, the stack of $r$-th roots of $\D$, denoted by $\X_r$, can be defined to be the moduli stack whose objects over a scheme $f: S\to \X$ are $$(M, \phi, \tau),$$
where 
\begin{enumerate}
\item
$M$ is a line bundle on $S$, 
\item
$\phi: M^{\otimes r}\simeq f^*\mathcal{O}_{\X}(\D)$ is an isomorphism,
\item
$\tau$ is a section of $M$ such that $\phi(\tau^r)$ is the tautological section of $f^*\mathcal{O}_{\X}(\D)$.
\end{enumerate}
The construction comes with a natural map $\X_r\to \X$, which is an isomorphism over $\X\setminus \D$. Let $$\D_r\subset \X_r$$ be the divisor lying over $\D\subset \X$.

In \cite{gs}, it is shown that the root construction above is essentially the only way stack structures can arise in codimension $1$. Naturally one can make the following
\begin{conjecture}[see \cite{ty1}, Conjecture 1.1]\label{conj:roots}
In the situation above, the Gromov-Witten theory of $\X_r$ is determined by the Gromov-Witten theory of $\X, \D$ and the restriction map $H^*_{CR}(\X)\to H^*_{CR}(\D)$. 
\end{conjecture}
However natural it may seem, Conjecture \ref{conj:roots} is not proven in the generality stated. The evidence known so far is the following 
\begin{theorem}[\cite{ty1}, Theorem 1.3]\label{thm:roots}
Conjecture \ref{conj:roots} holds true when $\D$ does not intersect the stacky locus of $\X$.
\end{theorem}
The above condition on $\D$ implies that $\D$ is a scheme. This is a rather restrictive assumption. However there is one instance in which this assumption holds even in repeated applications of root constructions: an orbifold curve $\C$ is obtained from its coarse moduli space, which is a nonsingular curve $C$, by a sequence of root constructions (see e.g. \cite{agv2}). Consequently, Theorem \ref{thm:roots} implies that the Gromov-Witten theory of an orbifold curve $\C$ is determined by the Gromov-Witten theory of its coarse moduli space, the nonsingular curve $C$. By the trilogy \cite{op1}, \cite{op2}, \cite{op3}, the Gromov-Witten theory of $C$ is completely determined. Therefore, the Gromov-Witten theory of an orbifold curve $\C$ is completely determined\footnote{Unfortunately, not explicitly determined at this moment.}. 

\subsection{Proof of a special case}  
We now give an outline of the arguments for Theorem \ref{thm:roots}. The first idea, which is quite standard in Gromov-Witten theory, is to use the deformation to the normal cone of $\D_r\subset \X_r$. This yields a degeneration of $\X_r$ to a union $$\X_r\cup_{\D_r}\Y$$ with $$\Y:=\mathbb{P}_{\D_r}(N_{\D_r/\X_r}\oplus \mathcal{O})$$ the projectivization of the normal bundle $N_{\D_r/\X_r}$ of $\D_r\subset \X_r$. The degeneration formula for orbifold Gromov-Witten theory, proven in \cite{af}, applies to express Gromov-Witten invariants of $\X_r$ in terms of relative Gromov-Witten invariants of the two pairs $$(\X_r, \D_r),\quad (\Y, \D_0\simeq \D_r),$$ where $\D_0\subset \Y$ is the $0$-section. Thus we need to determine Gromov-Witten invariants of the pairs $(\X_r, \D_r)$ and $(\Y, \D_r)$.

By \cite[Proposition 4.9]{af}, Gromov-Witten invariants of $(\X_r, \D_r)$ are determined by Gromov-Witten invariants of $(\X, \D)$. Because $\D$ is assumed to be a scheme, we can employ the proof of the ``relative-in-terms-of-absolute'' result in \cite{mp} to show that the relative Gromov-Witten theory of $(\X, \D)$ is determined by the Gromov-Witten theory of $\X, \D$ and the restriction map $H^*_{CR}(\X)\to H^*_{CR}(\D)$. 

It remains to treat the Gromov-Witten theory of $(\Y, \D_r)$. The torus $\mathbb{C}^*$ acts on $\Y$ by scaling the fibers. The fixed loci are $$\D_0, \D_\infty\subset \Y,$$ the ``$0$- and $\infty$-sections. Both divisors are isomorphic to $\D_r$. This $\mathbb{C}^*$ action on $\Y$ induces $\mathbb{C}^*$-actions on relevant moduli spaces, which allows us to apply virtual localization formula \cite{gp}, \cite{gv} to express Gromov-Witten invariants of $(\Y, \D_r)$ as integrals over $\mathbb{C}^*$-fixed loci. Analyzing the explicit virtual localization formula for relativce Gromov-Witten theory given in \cite{gv}, we find that the relative Gromov-Witten theory of $(\Y, \D_0)$ is determined by the Gromov-Witten theory of $\Y$ and the ``rubber theory\footnote{This is a variant of  relative Gromov-Witten theory.}'' for $(\Y, \D_0, \D_\infty)$.  

Virtual localization formula expresses Gromov-Witten invariants of $\Y$ in terms of integrals over the moduli space $\Mbar_{g,n}(\D_r, d)$ whose integrands involve the inverse $\mathbb{C}^*$-equivariant Euler class of certain $K$-theory class arising from the line bundle $N_{\D_r/\X_r}$. These quantities are examples of the so-called {\em twisted Gromov-Witten invariants} of $\D_r$, in the sense of \cite{tseng05}. The main theorem of \cite{tseng05}, usually called the orbifold quantum Riemann-Roch theorem, expresses these twisted Gromov-Witten invariants in terms of the usual Gromov-Witten invariants of $\D_r$ in a somewhat explicit fashion.

We are thus led to studying Gromov-Witten theory of $\D_r$. Since $\D_r$ is a root gerbe over $\D$ and $\D$ is a scheme, we can apply the results of \cite{ajt2, ajt3} (which are special cases of \cite{tt16} surveyed in Section \ref{sec:gerbe}) to explicitly write Gromov-Witten invariants of $\D_r$ in terms of Gromov-Witten invariants of $\D$.

It remains to study the rubber theory for $(\Y, \D_0, \D_\infty)$. Firstly, we apply the rigidification argument in \cite{mp} to write rubber invariants of $(\Y, \D_0, \D_\infty)$ as relativce Gromov-Witten invariants of $(\Y, \D_0, \D_\infty)$ with some ``target descendant'' insertions. We then apply an elaborated recursion called rubber calculus in \cite{mp} to remove these target descendants. This allows us to express rubber invariants of $(\Y, \D_0, \D_\infty)$ in terms of relative Gromov-Witten invariants of $(\Y, \D_0, \D_\infty)$.     

It can be seen from the construction that $\Y$ is obtained from a $\mathbb{P}^1$-bundle over $\D$, $\overline{\Y}\to \D$, by applying the $r$-th root constructions to its $0$- and $\infty$-sections $\overline{\D_0}, \overline{\D_\infty} \subset \overline{\Y}$. Therefore by \cite[Proposition 4.9]{af}, Gromov-Witten invariants of $(\Y, \D_0\cup \D_\infty)$ are determined by .Gromov-Witten invariants of $(\overline{\Y}, \overline{\D_0}\cup\overline{\D_\infty})$. Because $\D$ is a scheme, the quantum Leray-Hirsch theorem in \cite{mp} is applicable to the $\mathbb{P}^1$-bundle $\overline{\Y}\to \D$, implies that the Gromov-Witten theory of $(\overline{\Y}, \overline{\D_0}\cup\overline{\D_\infty})$ is determined by the Gromov-Witten theory of $\D$ and the restriction map $H^*_{CR}(\X)\to H^*_{CR}(\D)$.

The proof of Theorem \ref{thm:roots} is thus complete.

\subsection{Outlook}
Naturally, we would like to prove Conjecture \ref{conj:roots} in general. It is hopeful that the approach used in the proof of Theorem \ref{thm:roots} can be modified to address the general case, with no condition on $\D$. Indeed, it is not hard to see that several steps in the above outline do not use the condition on $\D$. However the condition on $\D$ is used in an essential way in the study of Gromov-Witten theory of $(\Y, \D_r)$. Therefore another approach must be taken.

In the author's opinion, the main difficulty lies in proving a quantum Leray-Hirsch result for Deligne-Mumford stacks of the form $\mathbb{P}_\D(L\oplus \mathcal{O})$ where $L$ is a line bundle over $\D$. The relative Gromov-Witten theory of $\mathbb{P}_\D(L\oplus \mathcal{O})$ relative to the $0$- and $\infty$-sections when the curve classes lie in the fiber of  $\mathbb{P}_\D(L\oplus \mathcal{O})\to \D$ requires solving the rubber theory for $1$-dimensional targets of the form $$([\mathbb{P}^1/G], [0/G]\cup [\infty/G]).$$ This is recently worked out in \cite{ty2}. However a quantum Leray-Hirsch result for $\mathbb{P}_\D(L\oplus \mathcal{O})$ relative to its $0$- and $\infty$-sections that covers curve classes of mixed type (i.e. combinations of fiber classes and classes from $\D$) is not yet proven.

It is also very desirable to have an {\em effective} version of Conjecture \ref{conj:roots}. More precisely, we would like to have some explicit formulas for Gromov-Witten invariants of $\X_r$ in terms of Gromov-Witten invariants of $\X$ and $\D$. It is reasonable to start with genus $0$. In various results on genus $0$ Gromov-Witten theory of a Deligne-Mumford stack $\X$, an essential role is played by the {\em J-function} of $\X$. This is a generating function of $1$-point genus $0$ descendant Gromov-Witten invariants of $\X$ which are important in the structures of genus $0$ Gromov-Witten theory. Thus it is very natural to ask the following 
\begin{question}
Is there an explicit formula for the $J$-function of $\X_r$, or a slice of the Lagrangian cones \cite{giv} defined by the genus $0$ Gromov-Witten theory of $\X_r$, given in terms of genus $0$ Gromov-Witten invariants of $\X$ and $\D$?
\end{question} 
To the author's knowledge, no answers to this Question is known in any examples. Since genus $0$ Gromov-Witten theory of toric stacks are explicitly calculated \cite{ccit}, it makes sense to try to study this Question when $\X$ is a toric stack and $\D$ is a toric prime divisor. Experiments in this toric setting may suggest possible answers to this Question in general.

\end{document}